\documentclass[11pt]{amsart}

\usepackage{amssymb,pstricks}
\begin{document}

\newtheorem{cor}{Corollary}
\newtheorem{de}[cor]{Definition}
\newtheorem{ob}[cor]{Remark}
\newtheorem{te}[cor]{Theorem}
\newtheorem{p}[cor]{Proposition}
\newtheorem{lema}[cor]{Lemma}
\newcommand{\dem}{{\noindent{\sc {Proof. }}}}

\newcommand{\Adelp}{{}^{\Phi}\!\!\cA_{\partial}}
\newcommand{\bX}{{\partial X}}
\newcommand{\cA}{\mathcal{A}}
\newcommand{\cC}{\mathcal{C}}
\newcommand{\cE}{\mathcal{E}}
\newcommand{\cF}{\mathcal{F}}
\newcommand{\ch}{\operatorname{ch}}
\newcommand{\cN}{\mathcal{N}}
\newcommand{\cR}{\mathcal{R}}
\newcommand{\cun}{\cC^{\infty}}
\newcommand{\cV}{\mathcal{V}}
\newcommand{\cz}{\mathbb{C}}
\newcommand{\Dp}{\Delta_{\Phi}}
\newcommand{\End}{\operatorname{End}}
\newcommand{\ffp}{{\rm ff}_{\Phi}}
\newcommand{\hol}{{\mathrm{hol}}}
\newcommand{\ind}{\operatorname{index}}
\newcommand{\nz}{\mathbb{N}}
\newcommand{\Ppm}{\Psi^{m_0}_{\Phi}(X)}
\newcommand{\Ppma}[1]{\Psi^{#1}_{\Phi}(X)}
\newcommand{\pnbX}{{}^{\Phi}N\bX}
\newcommand{\pnsbX}{{}^{\Phi}N^*\bX}
\newcommand{\Pnpm}{\Pnpma{m}}
\newcommand{\Pnpma}[1]{\Psi^{#1}_{{\rm sus}(\pnsY)-\varphi}(\bX)}
\newcommand{\pnY}{{}^{\Phi}NY}
\newcommand{\pnsY}{{}^{\Phi}N^*Y}
\newcommand{\pox}{{}^{\Phi}\Omega(X)}
\newcommand{\Ppa}[1]{\Psi_{\Phi}^{#1}}
\newcommand{\pssX}{{}^{\Phi}S^{*}X}
\newcommand{\ptsX}{{}^{\Phi}T^{*}X}
\newcommand{\ptX}{{}^{\Phi}TX}
\newcommand{\pXp}{{(\bX)^2_\varphi}}
\newcommand{\re}{\operatorname{Re}}
\newcommand{\res}{\operatorname{Res}}
\newcommand{\rhoN}{x}
\newcommand{\rpq}{\overline{\rz}_{+}}
\newcommand{\rz}{\mathbb{R}}
\newcommand{\tr}{\operatorname{tr}}
\newcommand{\Tr}{\operatorname{Tr}}
\newcommand{\Trdh}{\widehat{\Tr}_{\partial}}
\newcommand{\Trds}{\Tr_{\partial,\sigma}}
\newcommand{\Trsh}{\widehat{\Tr}_{\sigma}}
\newcommand{\vol}{{\mathrm{vol}}}
\newcommand{\vpX}{\cV_{\Phi}(X)}
\newcommand{\Xzp}{X^{2}_{\Phi}}
\newcommand{\Wres}{\operatorname{Wres}}
\newcommand{\zz}{\mathbb{Z}}

\newcommand{\mnote}[1]{\marginpar{\tiny\tt #1}}

\title{An index formula on manifolds with fibered cusp ends}
\author{Robert Lauter}
\author{Sergiu Moroianu}
\thanks{Moroianu was partially supported by a 
DFG-grant (436-RUM 17/7/01) and by 
the European Commission RTN HPRN-CT-1999-00118 \emph{Geometric Analysis}.}
\date{\today}
\begin{abstract}
We consider a compact manifold $X$ whose boundary is a locally trivial fiber 
bundle and an associated pseudodifferential algebra
that models fibered cusps at infinity. Using trace-like functionals that generate
the $0$-dimensional Hoch\-schild cohomology groups we express the index of a fully 
elliptic fibered cusp operator as the sum of a local 
contribution from the interior of $X$ 
and a term that comes from the boundary. This answers the index 
problem formulated in \cite{mame99}. We give a more precise answer in the 
case where the base of the boundary fiber bundle is $S^1$. In particular,
for Dirac operators associated to a product metric of the form
$g^X=\frac{dx^2}{x^4}+\frac{d\theta^2}{x^2}+g^F$ near the boundary $\{x=0\}$ 
with twisting bundle $T$ we obtain
$$\ind(A)=\int_X \hat{A}(X)\ch(T) -\frac{\lim_a\eta(A_{|\bX})}{2}$$
in terms of the integral of the Atiyah-Singer form
in the interior of $X$, and the adiabatic limit of the 
$\eta$-invariant of the restriction of the operator to the boundary.
\end{abstract}
\address{Fachbereich 17 - Mathematik\\ Universit\"at Mainz\\
D-55099 Mainz\\ Germany}
\email{lauter@mathematik.uni-mainz.de}
\address{Institutul de Matematic\u a al Academiei Rom\^ane\\ PO. Box 1-764\\
RO-70700, Bucharest\\ Romania}
\address{Universit\'e Paul Sabatier, UFR MIG,
118 route de Narbonne, 31062 Toulouse, France}
\email{moroianu@alum.mit.edu}
\maketitle

\section{Introduction}
\label{phi21}

Let $X$ be a compact manifold whose boundary is the total
space of a locally trivial fiber bundle
$\varphi:{\bX}\rightarrow Y$ of closed manifolds. 
Let $\rhoN:X\rightarrow\rpq$  be a defining function for
${\bX}$, i.e.\ ${\bX}=\{\rhoN=0\}$ and $d\rhoN$ does not vanish 
at ${\bX}$. One way to organize a (pseudo)differential analysis that 
reflects the geometry of $X$ is to choose a Lie algebra $\cV$ of vector fields 
on $X$, or more precisely a boundary fibration structure in the sense of Melrose
\cite{meicm}. 
The choice of the boundary fibration structure is by no means unique and 
different such structures on $X$ require in fact completely 
different analytic tools -- see, for instance, 
\cite{phiho1,de2,defr,maz91,mame99,mecom}.
In this paper we study the boundary fibration
structure determined by the Lie algebra $\vpX$ of
{\em fibered-cusp} or briefly {\em $\Phi$-vector fields} where a smooth 
vector field $V$ on $X$ belongs to $\vpX$ provided
$V$ is tangent to the fibers of $\varphi$ at the boundary and satisfies
$V\rhoN\in\rhoN^{2}\cun(X)$.
It is instructive to picture $\Phi$-vector fields in local coordinates.
Let $(x,y,z):X\supseteq U\rightarrow \rpq\times\rz^{n}_{y}\times\rz^{m}_{z}$ be local
product coordinates near the boundary such that $y, z$ are coordinates in $Y$,
respectively in the fiber $F$ of $\varphi$. Then $V\in\vpX$ can be written as
$$V|_{U}(x,y,z)=a(x,y,z)x^{2}\partial_{x}+
\sum_{j=1}^{n}b_{j}(x,y,z)x\partial_{y_{j}}+
\sum_{k=1}^{m}c_{k}(x,y,z)\partial_{z}$$
with coefficients $a$, $b_{j}$, and $c_{k}$ smooth down to $x=0$.
Taking these coefficients as local trivializations we see that the Lie
algebra $\vpX$ can be realized as the space of smooth sections of a 
smooth vector bundle $\ptX\rightarrow X$ that comes equipped with a natural
homomorphism $\ptX\rightarrow TX$.

The algebra
of $\Phi$-differential operators is by definition the enveloping algebra
${\rm Diff}^{*}_{\Phi}(X)$ of $\vpX$ over $\cun(X)$. 
A corresponding $\Phi$-pseudodifferential
calculus $\Ppma{*}$ has been constructed by Mazzeo and Melrose in \cite{mame99}.
First, $\Phi$-pseudodifferential operators act canonically 
on $\cun(X)$,
but since ${\rm Diff}^{*}_{\Phi}(X)$ as well as
$\Ppma{*}$ are $\cun(X)$-modules we can consider 
$\Phi$-(pseudo)differential operators acting between sections
of smooth vector bundles ${\mathcal E},{\mathcal F} \rightarrow X$ over $X$ and we write
${\rm Diff}^{*}_{\Phi}(X;{\mathcal E},{\mathcal F})$ resp.\
$\Psi^{*}_{\Phi}(X;{\mathcal E},{\mathcal F})$ for the corresponding spaces. 
Important examples of $\Phi$-differential operators are the Laplacian and the 
Dirac operators corresponding to \emph{exact $\Phi$-metrics}, for instance 
(\ref{gphi}) (see \cite{mame99}; 
these are certain complete metrics on the interior of $X$ which induce Euclidean
metrics on $\ptX$). 

As in the closed case, a $\Phi$-pseudodifferential operator
of order $m_0$ acts continuously as an operator of order $m_0$ 
on a scale of {\em $\Phi$-Sobolev spaces} $H_{\Phi}^{s}$, $s\in\rz$.
The $\Phi$-pseudodifferential operators $A$ that induce Fredholm operators
$H^{s}_{\Phi}\rightarrow H^{s-p}_{\Phi}$ have been characterized by
Mazzeo and Melrose \cite{mame99} as being those operators with invertible principal
symbol as well as invertible \emph{normal operator} (see Section \ref{phi22}). Such
operators are called fully elliptic. The index of a fully elliptic operator
$A$ is independent of the particular $s\in\rz$. A preliminary index formula 
for fully elliptic
operators has been obtained in \cite{phiho1} under the assumption that 
${\mathcal E}={\mathcal F}$ (this assumption is no restriction
when a fully elliptic operator exists; an isomorphism 
${\mathcal E}\to{\mathcal F}$ is given
by the principal symbol of the operator applied to a non-vanishing 
$\Phi$-vector field,
which exists whenever $\bX\neq\emptyset$).

Let us briefly recall the index formula.
We need several trace-like functionals on the $\Phi$-calculus whose 
definition has been adapted from a similar context in \cite{meni96c}.
Let $Q\in\Ppma{1}$ be a positive self-adjoint fully elliptic operator and 
$Q^{-\lambda}\in \Ppma{-\lambda}$ the family of complex powers constructed as in 
\cite{buc}.
Then for $A\in\Ppma{m_{0}}$, 
the operator $\rhoN^{z}AQ^{-\lambda}$ is of trace class
for $\re(z) > n+1$ and $\re(\lambda) >m_{0}+\dim(X)$, and the map
$(\lambda,z)\mapsto z\lambda\Tr \rhoN^{z}A Q^{-\lambda}$ admits a meromorphic 
extension to $\cz^{2}$ which is analytic near $(z,\lambda)=(0,0)$; thus,
we can define
\begin{eqnarray*}
\label{deftrace}
z\lambda \Tr \rhoN^{z}A Q^{-\lambda} &=:&
\Trds(A)+\lambda\Trdh(A)+ z\Trsh(A)\\
&&+\lambda^{2}W(\lambda,z)+
\lambda z W'(\lambda,z) +z^{2}W''(\lambda,z)\,,
\end{eqnarray*}
where the error terms $W,W',W''$ are holomorphic near $0\in\cz^{2}$.
The functionals obtained in this way have been studied in \cite{phiho1}.
 
\begin{te}[\cite{phiho1}]\label{iform}
Let $A\in\Psi_\Phi^{m_0}(M,\cE)$ be a fully elliptic $\Phi$-operator. Then
\begin{equation*}
\ind A  = \Trsh(A[B,\log Q]) - \Trdh([A,\log\rhoN]B)
\end{equation*}
where $B\in\Ppma{-p}$
is any inverse of $A$ up to trace class remainders, and
\begin{eqnarray*}
[B,\log Q]&:=&
\frac{d}{d\lambda} (Q^{-\lambda}BQ^{\lambda})_{|{\lambda=0}}
\in\Psi_\Phi^{-p}(X,\cE);\\
{}[\log\rhoN,A]&:=&\frac{d}{dz}(\rhoN^{z}A\rhoN^{-z})_{|{z=0}}\in
x\Psi_\Phi^{p-1}(X,\cE).
\end{eqnarray*}
\end{te}
Formally this is the same simplification of the computation of Melrose and Nistor
\cite{meni96c} as in \cite{de2, seronote}.
From \cite{phiho1} we know already that the first contribution to the index
is local, i.e.\ does not change if we modify $A$ by an operator of sufficiently
negative order, whereas the second contribution is global but
depends only on the behavior of the operators at the boundary, i.e.\
it does not change if we modify $A$ by an operator that vanishes to 
sufficiently high order at the boundary.
We identify in Proposition \ref{as}
the local term in Theorem \ref{iform} with the regularized Atiyah-Singer integral for
the index, defined in terms of heat kernel asymptotics. 

Our main interest lies in the case where the base $Y$ of the fiber bundle 
$\varphi:{\bX}\rightarrow Y$ is the circle $S^{1}$, for first-order 
differential operators 
$A$ modeled on Dirac operators. Choose a connection in the boundary fiber bundle, 
that is a rule for lifting the
horizontal vector field $\partial_\theta$. 
Choose a smooth metric on the interior of
$X$ that with respect to a product decomposition $\partial
X\times [0,1]\subset X$ close to the boundary 
looks like 
\begin{equation}\label{gphi}
g^X=\frac{dx^2}{x^4}+\frac{d\theta^2}{x^2}+g^F\,
\end{equation}
where $\theta$ is the variable in the circle and
$g^F$ is a family of metrics on the fibers. Such a metric is called a 
\emph{product $\Phi$-metric}; it induces an Euclidean metric on the vector 
bundle $\ptX$. Assume moreover that
${\mathcal E}_{|\bX}=E^+\oplus E^-$ and fix Hermitian metrics 
and connections $\nabla$ in $E^\pm$.

\begin{te}
\label{index}
Let $\cE, \cF\to X$ be Hermitian vector bundles and let
$$A:C^\infty(X,\cE)\to C^\infty(X,\cF)$$ 
be a first-order differential operator which in a product decomposition $\partial
X\times [0,1]\subset X$, $\cE_{|\bX}=E^+\oplus E^-$ looks like 
\begin{equation} \label{aexi}
A=\sigma\left((x^2\partial_x -\frac{x}{2})I_2+\left[\begin{array}{cc}
-ix\tilde{\nabla}_{\partial_\theta} &
D^*\\D&ix\tilde{\nabla}_{\partial_\theta}\end{array}\right]\right)
\end{equation}
where $D$ is a family of invertible operators on the fibers of $\varphi$, $D^*$ is the
formal adjoint of $D$,
$\sigma$ is an isometry $\cE_{|\bX}\to \cF_{|\bX}$, and 
$\tilde{\nabla}_{\partial_\theta}:=
\nabla_{\partial_\theta}+\frac{1}{4} \Tr(L_{\partial_\theta}g^F)$.
Then $A$ is Fredholm as an unbounded operator on
$L^2(X,\cE,\cF, g^X)$ and its index is given by
$$\ind (A)  = \overline{AS}(A)-\frac{\lim_a \eta(\delta_{x})}{2}.$$
Here $\overline{AS}(A)$ is the 
integral on $X$ of the pointwise 
supertrace of the heat kernel of $A$, which is a local expression in the full 
symbol of $A$ and in the metric, while $\lim_a \eta(\delta_{x})$ 
is the adiabatic limit (the limit as $x$ tends to
$0$) of the eta invariant of the ``boundary operator''
$$\delta_{x}:=\left[\begin{array}{cc}
-ix\tilde{\nabla}_{\partial_\theta} &
D^*\\D&ix\tilde{\nabla}_{\partial_\theta}\end{array}\right]:C^\infty(\bX,\cE)
\to C^\infty(\bX,\cE).$$
\end{te}
Note that we do not need to actually construct the heat kernel of $A$ 
in order to define the local index density. This density is 
slightly less than $L^1$, see Proposition \ref{intl1}. It is instructive to
compare this result with the classical Atiyah-Patodi-Singer index formula 
\cite{aps1}. 

\begin{cor}\label{cord}
Let $(X,g^X)$ be a compact spin manifold with boundary
and $T\to X$ a Hermitian vector bundle with constant metric and
connection near $\bX$. Then the twisted Dirac operator $A$
on $(X,g^X)$ is of the form \eqref{aexi}. Moreover, $A$ is Fredholm on 
$L^2(X,\cE,\cF,g^X)$ if and only if the 
family $D$ of Dirac operators on the fibers is invertible, and
$$\ind (A)  = \int_X \hat{A}(X)\ch(T) -\frac{\lim_a \eta(\delta_{x})}{2}.$$
\end{cor}
This follows immediately from the local index theorem (see \cite{bgv}) 
and Theorem \ref{index}. One checks directly in this case that the 
curvature of $g^X$ is a smooth $2$-form on $X$ with values in $\End(\ptX)$
so $\hat{A}(X)\ch(T)$ is a smooth form on $X$, thus in $L^1$ (see Proposition 
\ref{sm}).

\begin{cor}
Let $A$ be as in Theorem {\rm \ref{index}}. Then the integral of the index density
is an integer if and only if the determinant
bundle of the boundary family $D$ has trivial holonomy.
\end{cor}
This is a trivial application, since the adiabatic limit
of the eta invariant equals the logarithm of the holonomy (see Section 
\ref{eial} for the definitions). In particular, under the assumptions of 
Corollary \ref{cord} the Atiyah-Singer volume form
defines an integral cohomology class.

A related index formula has been obtained by Nye and Singer 
\cite{nysi} for the spin Dirac operator on $X=S^1\times\rz^3$ where the 
boundary fiber bundle is the projection $S^1 \times S^2\to S^2$. 
In the general case Vaillant
\cite{boristh} gave a formula for the index of the Dirac operator of a
$d$-metric on manifolds with
fibered cusps. Vaillant's formula contains the adiabatic 
limit of the eta invariant in the form computed by Bismut and Cheeger \cite{bc}.
It seems therefore likely that Theorem \ref{index} continues to hold for
boundary fiber bundles with higher dimensional base, as conjectured in \cite{nysi}.

The paper is structured as follows: in Section \ref{eial} we recall results 
about the eta invariant of self-adjoint operators. 
Section \ref{phi22} is devoted to an introduction to fibered cusp
pseudodifferential operators, with a focus on traces. 
The index theorem \ref{iform} is reviewed in Section \ref{aif}.
Finally the proof of Theorem \ref{index} 
occupies Section \ref{phi23}. A curious feature of the proof is the appearance 
and
then cancellation of the integral over $S^1$ of the determinant of $D^*D$ in the
boundary term of the index formula.

\vspace{0.5cm}
\noindent{\bf Acknowledgments:} 
The authors are indebted to Richard Melrose for drawing their attention to the
fascinating world of pseudodifferential algebras associated
to boundary fibration structures. His patience and generosity
have been of invaluable help for us. The first named author is grateful
to Michael Singer and Boris Vaillant
for some private lectures on the $\Phi$-calculus.

\section{Review of eta invariants and adiabatic limits}\label{eial}

The eta function of an elliptic self-adjoint differential 
operator $\delta$ was defined by Atiyah, Patodi and Singer \cite{aps1} as
$$\eta
(\delta,s):=\Tr\left((\delta^2)^{-\frac{s+1}{2}} \delta\right).$$
(assuming that $\delta$ is invertible).

Consider the family of operators on $\cun(\bX,\cE)$ 
defined for $0\leq x<\infty$ by
$$\delta_x:=\left[\begin{array}{cc}
-ix\tilde{\nabla}_{\partial_\theta} &
D^*\\D&ix\tilde{\nabla}_{\partial_\theta}\end{array}\right]$$
where $D$ is a family of elliptic invertible first-order differential 
operators on the fibers of $\bX$, acting from $E^+$ to $E^-$.
\begin{lema}\label{lem1.1}
The operator $\delta_x$ is symmetric on $L^2(\partial
X,\cE,\cE, \frac{d\theta^2}{x^2}+g^F)$.
\end{lema}

\dem We must show that 
$\tilde{\nabla}_{\partial_\theta}$ is skew-symmetric.
Let $s_1,s_2\in\cun(\bX,\cE)$. Then
\begin{eqnarray*}
\lefteqn{(\tilde{\nabla}_{\partial_\theta}s_1,s_2)
+(s_1,\tilde{\nabla}_{\partial_\theta}s_2)}&&\\&=&
\int_{\bX}\left((\tilde{\nabla}_{\partial_\theta}s_1,s_2)+
(s_1,\tilde{\nabla}_{\partial_\theta}s_2)\right)d g^F d\theta\\
&=&\int_{\bX}
(\partial_\theta(s_1,s_2)+\frac{(s_1,s_2)}{2}\Tr(L_{\partial_\theta}g))d g^F d\theta\\
&=&\int_{\bX} L_{\partial_\theta}((s_1,s_2)d g^F)d\theta
\\&=&\int_{S^1}\partial_\theta
\left(\int_{\bX/S^1}(s_1,s_2)d g^F\right)d\theta=0.
\end{eqnarray*}
\qed

We will see in Section \ref{phi23} that $\delta_x$ is also
invertible for small enough $x>0$.

The eta invariant of $\delta_x$ is by definition the regularized value of 
$\eta(\delta_x,s)$ at $s=0$. In fact, the eta function is regular at $s=0$
(see \cite{aps1}). By the adiabatic limit, denoted
$\lim_a\eta(\delta_{x})$, we mean the limit $\lim_{x\to 0} \eta(\delta_x).$
Intuitively it corresponds to separating the fibers of $\bX\to S^1$
in the limit since the Riemannian distance between distinct fibers tends to
infinity.

Recall the definition of the determinant line bundle of the family $D$ with the
Bismut-Freed connection. Since $D$ is invertible, $\det(D)$ 
is defined as the trivial bundle $\cz\times
S^1\rightarrow S^1$ with the connection $d+\omega^{BF}(0)$, where $\omega^{BF}(0)$
is the finite value at $s=0$ of the meromorphic family of $1$-forms
$$\omega^{BF}(s):=\Tr\left((D^*D)^{-\frac{s}{2}}D^{-1}
\tilde{\nabla}_{\partial_\theta}(D)\right).$$
Clearly, the holonomy of the Bismut-Freed connection is
$$\hol(\det(D),\omega^{BF}(0))=e^{-\int_{S^1}\omega^{BF}(0)}.$$
We define the logarithm of the holonomy of $\det(D)$ as
\begin{equation}\label{loghol}
\log\hol(\det(D))=-\int_{S^1}\omega^{BF}(0)\in i\rz.
\end{equation}
Determinant bundles and eta invariants are linked by the global anomaly formula
of Witten \cite{wit}, initially proved in \cite{bf2} for Dirac operators. 
The general result that we need is taken from \cite{ela,etalim}. 
\begin{te}\label{ldal}
The adiabatic limit of the eta invariant of $\delta_x$ satisfies
$${\lim}_a\eta(\delta_x)=-\frac{1}{i\pi}\log\hol(\det(D)).$$
\end{te}

\section{The fibered cusp calculus of pseudo-differential operators}
\label{phi22}

In this section we introduce the basic facts about the fibered cusp
calculus from \cite{phiho1} that are used
in the next sections. For a thorough treatment of the
fibered cusp calculus we refer the reader to \cite{mame99,boristh}.
We continue to use the notations from \cite{phiho1}. 

Blowing-up a submanifold $N$ of a smooth manifold $M$ means replacing 
$N$ by the set of its real normal directions inside $M$, i.e.\ the sphere bundle 
of its normal bundle; this procedure is equally defined for manifolds with 
corners. The result of the blow-up is a new manifold 
with corners of codimension possibly higher by $1$ than those of the 
initial manifold.

\subsection{The construction of $\Phi$-operators}
Let 
$$\Xzp:=[X\times X; \bX\times \bX, \pXp\times\{0\}]$$ 
be the fibered-cusp double space, obtained by an
iterated real blow-up from $X^2$ as follows: first blow up the corner 
$\bX\times \bX$ (at this stage we get the celebrated  $b$-double space 
$X^2_b$). The new boundary hyperface introduced by this blow-up is diffeomorphic
to $\bX\times\bX\times [-1,1]$ under the following map:
the class (modulo $\rz^*_+$) of the non-zero normal vector 
$(V_1,V_2)$ at $(y_1,y_2)$ maps to $\frac{V_1(x)+V_2(x)}{V_1(x)-V_2(x)}$.
Thus $\bX\times\bX\times\{0\}$ is a well-defined submanifold of 
$X^2_b$ provided we fixed the boundary-defining function $x$. 
The second stage of the construction  involves blowing-up the fiber diagonal 
$$\pXp:=\bX\times_\varphi\bX=\{(p,q)\in \bX\times\bX;\varphi(p)=\varphi(q)\}$$
of the boundary fiber bundle which by the discussion above is also a 
submanifold of $X^2_b$. The space $\Xzp$ comes equipped with a canonical smooth structure
as a manifold with corners of codimension at most $2$, and with a smooth 
blow-down map 
$$\beta:\Xzp\to X^2$$ 
which extends the identical diffeomorphism $({\Xzp})^\circ
=({X^2})^\circ$ of the interiors. The last face 
introduced by blow-up is called the $\Phi$-front face, denoted $\ffp$.
The lifted diagonal $\Dp$ is by definition the closure in $\Xzp$ of the preimage 
under $\beta$ of the interior of the diagonal in $X^2$.

The motivation of the construction is the fact \cite[Corollary 1]{mame99}
that the space of
$\Phi$-differential operators is canonically isomorphic to the space of 
distributions on $\Xzp$ supported on $\Dp$, conormal to $\Dp$ and 
extendable across $\ffp$, with
values in the bundle $\cF\boxtimes\cE^*\otimes \Omega'$. Here $\Omega'$ is the
pull-back through the projection on the right factor of the 
$\Phi$-density bundle $\Omega(\ptX)$. Note that 
$$\cun(X,\Omega(\ptX))=x^{-\dim(Y)-2}\cun(X,\Omega(TX)).$$
This singularity of order $\dim(Y)+2$
will play a great role in the rest of the paper. In fact, the main reason
for assuming $Y=S^1$ in Theorem \ref{index} is making this order of singularity 
small.
One defines then \cite{mame99} $\Ppma{m_0}$ as the space of linear operators 
$A:\dot{C}^\infty(X,\cE)\to\dot{C}^\infty(X,\cF)$ such that the lift 
$\kappa_A$ of the Schwartz kernel $k_A$ to $\Xzp$ is a classical conormal 
distribution of order $m_0$ on $\Xzp$ with values in 
$\cE^*\boxtimes\cF\otimes \Omega'$, vanishing 
rapidly to all boundary faces other than $\ffp$ and extendable across 
$\ffp$. These operators extend to bounded operators between
appropriate $\Phi$-Sobolev spaces. The fibered cusp calculus is closed under
composition \cite[Theorem 2]{mame99}.

If $\varphi$ is the identity map, fibered-cusp operators are nothing else than
\emph{scattering} operators \cite{mesca}. In that case the identifier $\Phi$
for double spaces, tangent bundles etc.\ will be replaced by $sc$.

\subsection{The normal operator}
There exist two symbol maps on $\Psi_\Phi(X)$, both multiplicative under composition of operators. 
One is the usual conormal principal symbol (living on $\ptsX$).
To describe the second symbol $\cN$, called the \emph{normal operator}, 
first assume that $Y=S^1$. 
In this case the interior of the $\Phi$-front face
$\ffp$ is the total space of a trivial $2$-dimensional real vector bundle over
$\pXp$. The two real directions correspond to the normal
direction to $\bX$ in $X$, and to the normal direction to 
$\pXp$ inside $\bX\times \bX$. 
Then $\cN$ is obtained by 
``freezing coefficients'' at $\ffp$
and then Fourier--transforming in the two real directions:
$$\cN(A):=\widehat{\kappa_A|_{\ffp}}.$$
Note that
in the general case there are $\dim Y+1$ suspending variables.
This new symbol map $\cN$ surjects onto the algebra $\Pnpma{\zz}$ 
of families of classical pseudo-differential operators on $\{Z_p\times
\rz^2\}_{p\in S^1}$ 
invariant with respect to translations in $\rz^2$
(\emph{$2$-suspended operators} in the terminology of \cite{meleta}) where $Z_p$ is the 
fiber over $p\in S^1$ of the boundary fiber bundle $\varphi:\bX\to S^1$.
An operator $A$ in $\Psi_\Phi(X)$ is called {\em elliptic} if
its principal conormal symbol is pointwise invertible, and {\em fully elliptic} if, in addition,
the corresponding normal operator $\cN(A)$ consists of a pointwise 
invertible family of pseudodifferential operators.

\subsection{The formal boundary symbol}
The principal symbol and the normal operator 
are invariantly defined. As for standard pseudodiffe\-rential operators
there exists a more refined notion of 
\emph{formal symbol} map, associating to an operator the Laurent series 
of its full symbol at the sphere at infinity inside the radial compactification
$\overline\ptsX$
(however, this symbol depends on choices except for its first term, 
the principal symbol).
Similarly, we associate to an operator its \emph {formal boundary symbol}
\begin{equation}
\label{TSC}
q:\Psi_\Phi(X) \rightarrow \Pnpma{\zz}[[x]]\,.
\end{equation}
To construct $q$, first choose a product decomposition
$\bX\times[0,\epsilon)\hookrightarrow X$ of $X$ near $\bX$ so that
$x(y,t)=t$. 
Let $X_\epsilon$ be the image of this map, and $Y_\epsilon:=Y\times[0,\epsilon)$.
Thus $X_\epsilon$ fibers over $Y_\epsilon$ via $\varphi\times Id$ with fiber 
type $F$ and $(X_\epsilon)^2_\varphi$ fibers over $(Y_\epsilon)^2_{sc}$ with 
fiber type $F\times F$.

Lift through $\beta$ the diagonal embedding 
$\pXp\times[0,\epsilon)\hookrightarrow
X^2$ to an embedding
\begin{equation}\label{efd}
\pXp\times[0,\epsilon)\hookrightarrow \Xzp.
\end{equation}
Note that $\pXp\times \{0\}$ maps identically to itself
as the zero section in $\ffp$. Moreover, the image of (\ref{efd})
is exactly the preimage of $\Delta_{sc}$ under the fibration 
\begin{equation}\label{fsc}
(X_\epsilon)^2_\Phi\to(Y_\epsilon)^2_{sc}.
\end{equation}
Thus the normal bundle to 
$\pXp\times[0,\epsilon)$ inside $\Xzp$ is the pull-back  via (\ref{fsc})
of the normal bundle to $\Delta_{sc}$,
which is canonically isomorphic to ${}^{sc}TY_\epsilon$. 
Consequently we use the notation $N(\pXp\times[0,\epsilon))={}^{sc}TX_\epsilon$.
The total space of ${}^{sc}TX_\epsilon|_{x=0}$
coincides with the interior of the front face.
Construct a collar neighborhood map $\mu:{}^{sc}TX_\epsilon\hookrightarrow \Xzp$
which on $\ffp$ is the identity.
Take a $\Phi$-operator $A$ and pull back its 
lifted Schwartz kernel $\kappa_A$ from $\Xzp$ to the total space of
${}^{sc}TX_\epsilon$ using $\mu$. Note that $\mu^*(\Omega')=
\Omega_{fiber,R}\otimes \pi^*(\Omega({}^{sc}TY_\epsilon))$
is the tensor product of the density bundle in the second factor
of $\pXp$ and the pull-back of the scattering density bundle from 
$Y_\epsilon$. At the 
zero section in ${}^{sc}TX_\epsilon$ we can identify the scattering density
factor with an Euclidean density on the fibers of ${}^{sc}TX_\epsilon$ 
(this density will alow us to take Fourier transforms in the fibers).
Identify $\cE,\cF$ and ${}^{sc}TX_\epsilon$ with their pull-backs from 
the zero-section of the front face, take the Taylor series of $\mu^*(\kappa_A)$
at $x=0$ (we can differentiate with respect to $x$ a distribution conormal to 
$\Delta_{\bX}\times [0,\epsilon)_x$) and 
then Fourier transform the coefficients in the fibers of
${}^{sc}TX_\epsilon$. This defines the map $q$.

\subsection{Product on formal boundary symbols}
The formal boundary map $q$ depends on choices of connections and 
trivializations except for its leading term which is just the normal operator. 
We denote by $*$ the product induced by $q$ on $\Pnpm[[x]]$.
Recall that in the standard pseudo-differential case 
we can choose the formal symbol
map so that the induced product on formal series of homogeneous symbols
(the so-called star product) takes the form
$$a(y,\xi)*b(y,\xi)=a(y,\xi)b(y,\xi)+\frac1i\sum
\partial_{\xi_i}a(y,\xi)\nabla_{\partial_{y_i}}b(y,\xi)+\ldots$$
Similarly in the context of Theorem \ref{index}
we can choose the formal boundary symbol map $q$ with the
properties:
\begin{eqnarray}
q(x\tilde{\nabla}_{\partial_\theta})&=&i\tau\label{q1}\\
q(x^2\partial_x)&=&i\xi\label{q2}\\
q(P)&=&P\nonumber
\end{eqnarray}
where $\tau$ is the suspending variable cotangent to the base $S^1$ 
of the boundary, $\xi$ is the suspending variable conormal 
to $\bX$ in $X$ and $P$ is a fiberwise differential operator 
which is constant in the fixed product decomposition near the boundary.
\begin{lema}
For $U, V\in\Pnpma{\zz}[[x]]$, 
the product induced by $q$ takes the form
\begin{equation}
U*V=UV+\frac{x}{i}\frac{\partial U}{\partial \xi}
\left(x\frac{\partial V}{\partial x}+\tau\frac{\partial V}{\partial \tau}
\right)+\frac{x}{i}\frac{\partial U}{\partial
\tau}\tilde{\nabla}_{\partial_\theta}(V) +O(x^2).\label{prodsus}
\end{equation}
where the product in the right-hand side is the standard product
of power series with coefficients in the algebra $\Pnpma{\zz}$.
\end{lema}
\dem It is enough to prove the formula for $\Phi$-differential operators 
since the
product is given by bi-differential operators with polynomial coefficients
(see e.g.~\cite[Proposition 3.11]{defr} for details of this argument in a 
similar context). 
But ${\rm Diff}^{*}_{\Phi}(X,\cE)$ is generated near $\bX$ by
$ix\tilde{\nabla}_{\partial_\theta}$, $ix^2\partial_x$ and by
differential operators
along the fibers. The lemma follows easily from (\ref{q1}) and (\ref{q2}) 
since it is valid on the generators.
\qed

\subsection{Traces densities of $\Phi$-operators}
Of main interest for us are traces of $\Phi$-operators. We study them using the 
more refined notion of trace density. It is a standard fact that on a 
closed manifold $M$, any operator $A\in\Psi^{\lambda}(M,\cE)$ with 
$\re(\lambda)<-\dim(M)$ is of trace class, and
$$\Tr(A)=\int_\Delta \tr(k_A|_\Delta)$$
(the Schwartz kernel $k_A$ is continuous on $M\times M$ and its restriction
to the diagonal is a smooth $1$-density with values in $\End(\cE)$). 
The same remains true for 
$\Phi$-operators modulo an integrability issue at the boundary. By pulling-back 
via the blow-down map we write
$$\Tr(A)=\int_{\Dp}\tr(\kappa_A|_{\Dp}).$$
We identify $\Dp$ with $X$ via the right projection.
The restriction to $\Dp$ of the density bundle $\Omega'$ is precisely
$\Omega(\ptX)=x^{-\dim(Y)-2}\Omega(X)$.
It is then clear that $A$ is of trace class if and only if 
$A\in x^z\Psi_\Phi^\lambda(X,\cE)$ with $\re(z)>\dim(Y)+1, \re(\lambda)<-\dim(X)$.
The density $\tr(\kappa_A|_{\Dp})$, viewed as a density on $X$, is called the
\emph{trace density} of $A$.

Recall that for the definition of the formal boundary symbol we have chosen a 
product decomposition of $X$ near $\bX$, as well as a local 
isomorphism of $\cE$ with its pull-back from $\bX$. We can 
thus expand $\kappa_A|_{\Dp}$ in powers of $x$ near $x=0$.
Let $\pnbX$ be the restriction of the vector bundle 
${}^{sc}TX_\epsilon$ to $\bX\times\{0\}\subset \ffp$ and $\pnsbX$ its dual.

\begin{p}
Let $A\in x^z\Psi_\Phi^\lambda(X,\cE)$ with $\re(\lambda)<-\dim(X)$. Then
$$\kappa_A|_{\Dp}\sim_{x\to 0} \frac{1}{(2\pi)^{\dim(Y)+1}}
\int_{\pnsbX/\bX} q(A) 
\frac{\omega_{sc}^{\dim(Y)+1}}{(\dim(Y)+1)!}.$$
\end{p}
\dem This is precisely the Fourier inversion formula in each fiber
of $\pnbX$, where $\omega_{sc}$ is the canonical (singular) symplectic form on
${}^{sc}T^*Y_\epsilon$ pulled back to $\pnsbX$. 
\qed

In the case $Y=S^1$ this takes a somewhat simpler form.
\begin{cor}\label{kdq}
Let $A\in x^z\Psi_\Phi^\lambda(X,\cE)$ with 
$\re(\lambda)<-\dim(X)$ and assume $Y=S^1$. Then 
$$\kappa_A|_{\Dp}\sim_{x\to 0} \left(\frac{1}{(2\pi)^{2}}
\int_{\rz^2} q(A)(\tau,\xi)d\tau d\xi\right)\frac{d\theta dx}{x^3}.$$
\end{cor}
\dem Notice that $\pnbX$ is a trivial $\rz^2$ bundle in this case, while
$\omega_{sc}=\frac{d\xi\wedge dx}{x^2}+\frac{d\tau\wedge d\theta}{x}$.
\qed

If $\cz\ni\lambda\mapsto A(\lambda)\in\Psi_\Phi^\lambda(X,\cE)$ is an entire family of 
$\Phi$-operators then both $\kappa_A|_{\Dp}$ and $ q(A)|_\pnsbX$ 
extend meromorphically to $\cz$ with at most simple poles. Moreover,
Corollary \ref{kdq} remains valid for the meromorphic extensions 
for all values of $\lambda$ (at a pole, the expansions are valid 
for the residues and for the regularized parts separately).

We close this section with a description in terms of the map $q$ 
of the trace functional $\Trdh$ defined in the Introduction 
(compare with \cite[Proposition 7.6]{phiho1}). 

\begin{lema}\label{lem4}
Let $A\in\Psi_\Phi^{m_0}(X)$.
Then $\Trdh(A)$ is explictly given by
$$\Trdh(A)=\frac{1}{(2\pi)^2}\left(\int_{S^1\times\rz^2}
\Tr(q(AQ^{-\lambda})_{[-2]})d\theta d\tau d\xi\right)_{\lambda=0}$$
where $\Tr$ in the right-hand side denotes the trace of operators on the 
fibers of $\varphi:\bX\to S^1$, $(\cdot)_{[k]}$ is the coefficient of 
$x^{-k}$ in a power expansion in $x$,
and $(\cdot)_{\lambda=0}$ stands for the regularized value at $\lambda=0$.
\end{lema}
\dem We can assume that $\kappa_{AQ^{-\lambda}}$ is supported in 
$(X_\epsilon)^2_\Phi$. Then 
\begin {eqnarray*}
\Trdh(A)&=&\res_{z=0}\Tr(x^zAQ^{-\lambda})_{|\lambda=0}\\
&=&\res_{z=0}\int_{\Dp} \tr(\kappa_{x^zAQ^{-\lambda}})_{|\lambda=0}\\
&=&\res_{z=0}\int_0^\epsilon \int_{\bX} 
x^z \tr(\kappa_{AQ^{-\lambda}})_{|\lambda=0}
\end{eqnarray*}
Obviously, 
$$\res_{z=0}\left(\int_0^\epsilon x^{z-k}dx\right)$$
equals $1$ if $k=1$ and 
$0$ otherwise, so
\begin {eqnarray*}
\Trdh(A)&=&\left(\int_{\bX} 
\tr(\kappa_{AQ^{-\lambda}})_{[1]}\right)_{|\lambda=0}\\
&=&\left( \int_{S^1}\int_{\bX/S^1} 
\left(\frac{1}{(2\pi)^{2}}
\int_{\rz^2}\tr q(AQ^{-\lambda})_{[-2]}d\tau d\xi\right)d\theta\right)_{|\lambda=0}
\end{eqnarray*}
(we used Corollary \ref{kdq} and the remark following it 
in the last equality). Now 
$$\int_{\bX/S^1} \tr q(AQ^{-\lambda})_{[-2]}=\Tr( q(AQ^{-\lambda})_{[-2]})$$ 
so the result follows by Fubini's theorem. 
\qed

\section{The abstract index formula}
\label{aif} 

The results of this section hold for general boundary fiber bundles, i.e.~not
necessarily with base $S^1$.

Fibered-cusp operators have two types of principal symbols. Accordingly, 
elliptic regularity has a
new aspect in fibered-cusp theory concerning regularity at the boundary.
The following lemma and its proof are quite standard; we include them for 
future reference.

\begin{lema} \label{fef}
Let $A\in\Psi_\Phi(X,\cE,\cF)$ be fully elliptic. Then the 
$L^2$ solutions of $A\psi=0$ belong to $x^\infty 
C^\infty(M,{\mathcal E})$.
\end{lema}
\dem Since $A$ is fully elliptic there exists 
a parametrix $B$ of $A$ inverting $A$ up to
$R\in x^\infty \Psi_\Phi^{-\infty}(X,{\mathcal E})$.
Let $\psi$ be a distributional solution of the pseudo-differential equation
$A\psi=0$. It follows
$$0=BA\psi=(I+R)\psi=\psi+R\psi$$
so $\psi=-R\psi$. But $R\in x^\infty \Psi_\Phi^{-\infty}(X,{\mathcal E})$  
implies $R\psi\in x^\infty 
C^\infty(M,{\mathcal E})$.
\qed

Since $R$ is compact on $L^2_\Phi$ (the compact operators in
$\Ppma{\zz}$ are precisely those in $x\Ppma{-1}$) it follows that $\ker A$ is finite 
dimensional and moreover the orthogonal projection
$P_{\ker A}$ belongs to the ideal
$x^\infty \Psi_\Phi^{-\infty}(X,{\mathcal E})$.
We can define therefore invertible $\Phi$-operators 
\begin{eqnarray*}
Q_1&:=&(AA^*+P_{\ker A^*})^{1/2};\\
Q_2&:=&(A^*A+P_{\ker A})^{1/2}.
\end{eqnarray*}
Note that $q(Q_1)=q(AA^*)^{1/2}$. Let 
$$B:=A^*Q_1^{-1}$$
be a parametrix of $A$. 

\subsection{The index formula}
Let us reprove the index formula from \cite{phiho1}. Assume for simplicity 
that $A$ is of order $1$. For technical reasons we would like to work with 
operators acting from $\cE$ to itself.
Recall that $\bX\neq\emptyset$ implies the existence of a non-vanishing
vector field on $X$ (the obstruction to the existence of such a vector field 
lives in $H^{\dim(X)}$ and this space is $0$ when $\bX\neq\emptyset$). There 
exist non-canonical isomorphisms between $TX$, $\ptX$ and $\ptsX$,
and thus a non-vanishing section in $\ptsX$. The principal conormal symbol of $A$ 
evaluated on this section
gives an isomorphism $u$ between $\cE$ and $\cF$. Finally, $v:=u^*(uu^*)^{-1/2}$
is an isometry $\cF\to\cE$. Thus $U:=vA\in\Psi_\Phi^1(X,\cE)$ has the property
\begin{eqnarray*}
&&R_1:=(UU^*+P_{\ker U^*})^{1/2}
=vQ_1 v^*\\
&&R_2:=(U^*U+P_{\ker U})^{1/2}=Q_2.
\end{eqnarray*}
Set $V:=Bv^*$. Note the commutations $UR_2^{-\lambda}=R_1^{-\lambda}U$, 
$VR_1^{-\lambda}=R_2^{-\lambda}V$. The index formula is obtained as follows:
\begin{eqnarray}
\ind(A)&=&\ind(U)\nonumber\\
&=&\Tr(UV-VU)\nonumber\\
&=&\Tr(x^z(UV-VU)R_2^{-\lambda})_{\lambda=0,z=0}\nonumber\\
&\stackrel{(\ref{fint}.1)}{=}&\Tr(x^zUVR_2^{-\lambda}-Ux^z V R_1^{-\lambda})_{\lambda=0,z=0}\nonumber\\
&=&\Tr(x^z UV(R_2^{-\lambda}-R_1^{-\lambda})+[x^z,U]VR_1^{-\lambda})_{\lambda=0,z=0}
\nonumber\\
&\stackrel{(\ref{fint}.2)}{=}&(\Tr(x^z 
(R_2^{-\lambda}-R_1^{-\lambda}))_{\lambda=0,z=0}\nonumber\\&&
+\Tr([x^z,A]B
Q_1^{-\lambda})_{\lambda=0,z=0}.\label{fint}
\end{eqnarray}
For (\ref{fint}.1) we have used the fact that $\Tr[C,D]=0$ for $C\in x^c\Ppma{a}$,
$D\in x^d \Ppma{b}$ with $a,b,c,d\in\cz$, $\re(a+b)<-\dim(X)$,
$\re(c+d)>1+\dim(Y)$.
Thus (\ref{fint}.1) is true for large real parts of
$\lambda,z$, hence for any $\lambda,z\in\cz$ by unique continuation. 
(\ref{fint}.2) holds at $\lambda=0,z=0$ because $UV-1\in x^\infty \Ppma{-\infty}$
while $v,v^*$ cancel in the trace. 

\subsection{The interior term}
We claim that the first term in (\ref{fint}) is the regularized integral on 
$X$ of a local expression in the full symbol expansion of $U$. Indeed, for 
$j=1,2$ let $r_j(\lambda)\in C^\infty(X,\End(\cE)\otimes\pox)$ be the 
meromorphic extension of the lifted Schwartz kernel $\kappa_{R_j^{-\lambda}}$ 
restricted to $\Dp$. As in the case of
closed manifolds it is easy to see
that $r_j(\lambda)$ is regular at $\lambda=0$. The first term
in (\ref{fint}) is 
${\int_X x^z\tr(r_2(0)-r_1(0))}|_{z=0}.$ 
By Corollary \ref{kdq} and the remark after it, the density $\tr(r_2(0))$
has a Laurent expansion at $x=0$ starting with $x^{-3}$ so the previous
integral is absolutely integrable and 
holomorphic in $z$ for $\re(z)>2$ and extends to 
$\cz$ with possible simple poles at $z=2-\nz$.
Now 
\begin{eqnarray*}R_2^{-\lambda}-R_1^{-\lambda}&=&
U[R_2^{-\lambda},V]+O(\lambda) x^\infty\Ppma{-\infty}\\
&=&\lambda U[V,\log R_2]R_2^{-\lambda}+O(\lambda)x^\infty\Ppma{-\infty}+O(\lambda^2)
\end{eqnarray*}
where $O(\lambda)$ denotes an analytic
multiple of $\lambda$ near $\lambda=0$. Clearly then 
\begin{equation}\label{tl}
(\Tr(x^z (R_2^{-\lambda}-R_1^{-\lambda}))_{\lambda=0,z=0}= \Trsh(U[V,\log R_2]).
\end{equation}
On the more refined level of trace densities, by \cite[Proposition 7.4]{phiho1},
$$r_2(0)-r_1(0)=\frac{1}{(2\pi)^N}\int_{\pssX/X}\sigma_{[-N]}(U[V,\log R_2])
\imath_\cR \omega_\Phi^N$$ 
is given in terms of the component of homogeneity $-\dim(X)$ of the 
formal symbol of $U[V,\log R_2]$, so clearly depends only on the 
jets of the full symbol of $U$. 

\subsection{The boundary term}
Similarly for the second term from (\ref{fint})
we have 
$$[x^z,A]B=zx^z[\log x,A]B+O(z^2)$$
and $\Tr(O(z^2))=O(z)$ so
\begin{equation}\label{tb}
\Tr([x^z,A]BQ_1^{-\lambda})_{\lambda=0,z=0}=-\Trdh([A,\log x]B).
\end{equation}
By Lemma \ref{lem4} this last quantity is concentrated at the boundary.

So far, combining (\ref{fint}), (\ref{tl}) and (\ref{tb}) 
we have proved the general index Theorem \ref{iform}. Note that in 
(\ref{tl}) we can assume $v=1$ (and so $U=A$, $V=B$, $Q=R_2$) since in Theorem \ref{iform} we suppose
$\cE=\cF$. 

\subsection{Relationship with heat kernel expansions}
Here is a more familiar interpretation of the local term (\ref{tl}). 
It is worth stressing that we do not prove a heat kernel expansion for $A^*A$.
Rather we use the existence of heat kernel expansions for pseudo-differential
operators on closed manifolds as well
as the locality of the two quantities we want to relate. 
\begin{p} \label{as}
The local quantity $\tr(r_2(0)-r_1(0))$ equals the index density, defined as the 
universal expression in the 
jets of the full symbol of $U$ which gives the pointwise supertrace of the 
constant term in the heat kernel expansion.
\end{p}
\dem Fix a point $p$ in the
interior of $X$ and modify the operator $U$ far from $p$ so that it extends to an
elliptic operator on the double of $X$. Denote the extensions to $2X$ by the same
letters as before. Then for $j=1,2$ use the Mellin 
transformation formula
$$\Gamma(\lambda)R_j^{-2\lambda}=\int_0^\infty t^{\lambda-1}e^{-tR_j^2}dt$$
to identify the value at $\lambda=0$ of the analytic extension of the Schwartz kernel 
of $R_j^{-\lambda}$ on the diagonal with the coefficient of $t^0$ in the asymptotic
expansion as $t\searrow 0$ of the Schwartz kernel 
of $e^{-tR_j^2}$ on the diagonal.
Remember that $R_2=Q_2$, and observe that $Q_1^2$ and $R_1^2$ are 
conjugate via $v$ so the pointwise trace of their Schwartz kernels is the same.
\qed

Such a formula for the local term is not surprising in index theory. Our point
is getting it without having to construct heat kernels for $\Phi$-operators. 
In this respect the approach via complex powers,
which are already objects in the calculus, presents a great advantage.

Although the local term (\ref{tl}) is smooth on $X$ up to the boundary as 
a $\Phi$-density, its integral might in principle diverge since as an 
usual density on $X$ it has a singularity of order $3$ at $x=0$. Thus,
we cannot set directly $z=0$ in the above evaluations. To prove Theorem 
\ref{index} we must identify the boundary term with the adiabatic limit of the 
eta invariant and show that the index density is integrable 
in a restricted sense. We will do this in Section \ref{phi23}. 

However, we can prove directly that when $A$ is a twisted Dirac
operator corresponding to the metric $g^X$, and with twisting bundle $T$
with metric and connection constant in $x$ in a neighborhood of the boundary,
the local term in the index formula equals the integral of the
Atiyah-Singer density $\hat{A}(X,g^X)\ch(T)$ \emph{without regularization}.
Indeed,

\begin{p} \label{sm}
The Riemannian curvature $R$ of $(X,g^X)$ induces a smooth $2$-form on 
$X$ with values in $\End(\ptX)$ down to $x=0$. Thus 
$\hat{A}(X,g^X)$ and $\ch(T)$ are smooth forms on $X$.
\end{p}
\dem Let $e_1,\ldots,e_m$ be a local orthonormal frame in the fibers of $\phi$.
It is straightforward to compute $R(g^X)$ using the local orthonormal
$\Phi$-vector fields $x^2\partial_x$, $x\partial_\theta$ (lifted to $\bX$ using 
the connection involved in the definition of $g^X$),
$e_1,\ldots,e_m$. We observe that while for instance
$R(\partial_x,e_i)\partial_\theta$
diverges like $1/x$ as $x\to 0$, the induced action of $R$ on $\ptX$ is smooth 
down to $x=0$. To conclude that $\hat{A}(X,g^X)$ is smooth
it is enough to prove that $\tr(R^k)$ is smooth down to $x=0$
for all $k\in\nz$. Of course it does not 
matter for the trace if we view the $2$-form $R$ as acting on $TX$ or on
$\ptX$, so the conclusion follows. $\ch(T)$ is obviously smooth in $x$ 
(it is in fact constant in $x$ in a neighborhood of $x=0$).
\qed

\section{The index of first-order $\Phi$-differential operators}
\label{phi23} 

For the rest of the paper we assume that $Y=S^1$ and that $A$ satisfies the hypothesis
of Theorem \ref{index}. Note that $x^2\partial_x-x/2$ is skew-symmetric with respect
to the metric $g^X$. 
It follows directly from the properties of the quantization $q$ that
\begin{eqnarray}
q(\delta_x^2)&=&\left[\begin{array}{cc}
\tau^2+D^*D& -ix\tilde{\nabla}_{\partial_\theta}(D^*)\\
ix\tilde{\nabla}_{\partial_\theta}(D)&\tau^2+DD^* \end{array}\right]\nonumber\\
q(A)&=&\sigma\left[\begin{array}{cc}
i\xi-\frac{x}{2}+\tau & D^*\\D&i\xi-\frac{x}{2}-\tau\end{array}\right]\nonumber\\
q(A^*)&=&\left[\begin{array}{cc}
-i\xi+\frac{x}{2}+\tau &
D^*\\D&-i\xi+\frac{x}{2}-\tau\end{array}\right]\sigma^*\label{qas}\\
\cN(AA^*)&=&\sigma\left(\left[\begin{array}{cc}
\xi^2+\tau^2+D^*D& 0\\0&\xi^2+\tau^2+DD^*
\end{array}\right]\right)\sigma^*\label{cnaa}\\
q(AA^*)&=&\sigma\left(\left[\begin{array}{cc}
\xi^2+\tau^2+D^*D& 0\\0&\xi^2+\tau^2+DD^* \end{array}\right]\right.\nonumber
\\&&+\left. x\left[\begin{array}{cc}
i\xi+\tau&-i\tilde{\nabla}_{\partial_\theta}(D^*)\\
i\tilde{\nabla}_{\partial_\theta}(D)&
i\xi-\tau \end{array}\right]+\frac{x^2}{4}\right)\sigma^*.\label{aas}
\end{eqnarray}
Set $\Delta_+:=D^*D$, $\Delta_-:=DD^*$, $\Delta:=\left[\begin{array}{cc}
\Delta_+&0\\0&\Delta_-\end{array}\right]$.

\begin{lema}\label{lemma4.1}
The operator $A$ defined by \eqref{aexi} is a fully elliptic $\Phi$-operator if and
only if the family $D$ is invertible.
\end{lema}

\dem It is clear from (\ref{cnaa}) that 
$\cN(AA^*)=\sigma(\xi^2+\tau^2+\Delta)\sigma^*$ is elliptic 
and non-negative as a $2$-suspended operator; moreover it is invertible 
for each 
value of the parameters $\tau,\xi\in\rz$ (thus invertible as a 
suspended operator, see \cite{meleta}) if and only if $D$ is invertible.
\qed

\begin{lema}\label{lemma4.2} \cite{ela, etalim}
The differential operator
$\delta_x\in \mathrm{Diff}^1(\bX,{\mathcal E})$ is 
invertible for $0<x<\epsilon$ for some $\epsilon>0$.
\end{lema}

\dem (sketch) We can view $\delta_x$ as an adiabatic family of operators (i.e.\ an adiabatic
differential operator in the sense of \cite{etalim}). The adiabatic 
normal operator of this family is invertible as in Lemma \ref{lemma4.1}, 
so exactly like in Lemma \ref{fef} there exists an inverse 
$\mu_x\in \Psi_a^{-1}(\bX, \cE)$ modulo $x^\infty\Psi_a^{-\infty}(\bX,\cE)$:
$$\delta_x \mu_x=I+r_x$$
Now the residual adiabatic ideal $x^\infty\Psi_a^{-\infty}(\bX,\cE)$
equals the space of rapidly vanishing families
of smoothing operators on $\bX$ as $x\to 0$. Thus 
$r_x\to 0$ inside bounded operators as $x\to 0$. The conclusion follows
for $\epsilon$ chosen small enough so that $\|r_x\|<1, \forall x<\epsilon.$
\qed

\begin{p}\label{prop8.4}
The boundary term 
$\Trdh([A,\log \rhoN]B)$ from the index formula \eqref{fint}
equals $-\frac{i}{2\pi}\log(\hol(\det D))$.
\end{p}
\dem  
From (\ref{aexi}) we know $q([A,\log \rhoN])=x\sigma$.  We claim that 
we can assume $\sigma=1$. Let $S,T\in\cun(X,\End(\cE))\subset\Psi_\Phi^{0}(M,\cE)$ be such that
$q(S)=\sigma, q(T)=\sigma^{-1}$. Since $q(ST)=1$ it follows 
\begin{eqnarray*}
\Tr(x^z[A,\log \rhoN]BQ_1^{-\lambda})
&=&\Tr(x^z ST[A,\log \rhoN]BQ_1^{-\lambda})
+O(z^0)\\
&=&\Tr(x^z T[A,\log \rhoN]BQ_1^{-\lambda}S)+O(z^0).
\end{eqnarray*} 
Observe that $q(T[A,\log \rhoN])$ and $q(BQ_1^{-\lambda}S)$ do not contain 
$\sigma$ anymore; on the other hand the term regular in $z$ at $z=0$ 
does not affect the residue, which proves our claim.
From Lemma \ref{lem4}, 
\begin{eqnarray}
\lefteqn{\Trdh([A,\log \rhoN]B)}&&\nonumber\\
&=&\frac{1}{(2\pi)^2}\int_{S^1\times\rz^2}
\Tr\left( q([A,\log \rhoN]BQ_1^{-\lambda})_{[-2]} \right) d\tau d\xi 
d\theta_{|\lambda=0}\nonumber\\
&=&\frac{1}{(2\pi)^2}\int_{S^1\times\rz^2}
\Tr\left(q(A^*(AA^*)^{-\frac{\lambda}{2}-1})_{[-1]} \right) d\tau
 d\xi d\theta_{|\lambda=0}
\label{etlm}
\end{eqnarray}
where we use the formulas (\ref{qas}), (\ref{aas}) 
for $q(A^*), q(AA^*)$ with $\sigma=1$.
There are three types of terms 
occurring in (\ref{etlm}) as explained below and we write accordingly
$$\Trdh([A,\log \rhoN]B)=(I(\lambda)+II(\lambda)+III(\lambda))_{|\lambda=0}.$$ 

\subsection{The terms of type $I$}
First there are those terms where $q(A^*)$ and
$\cN((AA^*)^{-\lambda/2-1})$
are composed according to the product rule (\ref{prodsus}).
Since $\frac{\partial q(A^*)}{\partial \xi}$ and $\frac{\partial
q(A^*)}{\partial \tau}$ are constants, the corresponding integrands are exact forms
so the terms containing them vanish.
The non-vanishing terms of type $I$ come from
$\frac{x}{2}\cN((AA^*)^{-\lambda/2-1})$ and $\frac{x\tau}{i}
\frac{\partial \cN(A^*)}{\partial \xi} \frac{\partial
\cN(AA^*)}{\partial \tau}$.
Using (\ref{cnaa}) and polar coordinates in the $(\tau,\xi)$ plane we get 
\begin{eqnarray}
I_1(\lambda)&=&\frac{1}{(2\pi)^2}\int_{S^1\times\rz^2}
\frac{1}{2}\Tr\left(\cN(AA^*)^{-\lambda/2-1}\right)d\theta d\tau d\xi\nonumber\\
&=&\frac{1}{4\pi\lambda} \int_{S^1}\Tr\left(\Delta^{-\frac{\lambda}{2}}\right)d\theta.
\label{u1}
\end{eqnarray}
Thus at $\lambda=0$ we get the average of the logarithm of the determinant
of the family $D$. Although in the end this term will cancel away, it is 
worth recalling the definition of the
determinant of $\Delta$, not to confuse with the determinant 
line bundle with connection 
$(\det D, d+\omega^{BF})$ defined in Section \ref{eial}. 
The zeta function of any positive pseudo-differential
operator is regular at $\lambda =0$;  the logarithm of the determinant of $\Delta$
is defined as the derivative $\zeta'(\Delta,0)$. This derivative
clearly equals the finite part at $\lambda=0$ of 
$-\frac{1}{\lambda}\Tr(\Delta^{-\frac{\lambda}{2}})$. 

Similarly we get
\begin{eqnarray}
I_2(\lambda)&=&\frac{1}{(2\pi)^2}\int_{S^1\times\rz^2}\Tr\left(
\frac{\tau}{i}\frac{\partial q(A^*)}{\partial\xi}
\frac{\partial \cN(AA^*)^{-\frac{\lambda}{2}-1}}{\partial \tau}
\right)d\theta d\tau d\xi \nonumber\\
&\stackrel{(\ref{ti}.1)}{=}&\frac{1}{(2\pi)^2}\int_{S^1\times\rz^2}
\Tr(\tau^2+\xi^2+\Delta)^{-\frac{\lambda}{2}-1}
d\theta d\tau d\xi\nonumber\\
&\stackrel{(\ref{ti}.2)}{=}&\frac{1}{2\pi\lambda} \int_{S^1}\Tr\left(\Delta^{-\frac{\lambda}{2}}\right)d\theta.
\label{ti}
\end{eqnarray}
(in (\ref{ti}.1) we used the formulas (\ref{qas}) and (\ref{cnaa}) with $\sigma =1$
while in (\ref{ti}.1) we used polar coordinates in the $(\tau, \xi)$ plane).

\subsection{The terms of type $II$}
The second type of terms in (\ref{etlm})
come from the coefficient of $x$ in $(\xi^2+\tau^2+\Delta)^{-\lambda/2-1}_{[-1]}$, 
where the power is taken with respect to the product
(\ref{prodsus}). This diagonal matrix is not explicitly computable, however the 
trace of $\cN(A^*)$ times it \emph{is}, because of two facts:
\begin{itemize}
\item The diagonal of $\cN(A^*)$ is made of central elements
modulo $x$.
\item The partial derivatives of $\xi^2+\tau^2+\Delta$
with respect to $\xi$ and $\tau$, are central
elements in $\Ppma{\zz}/x^\infty \Ppma{\zz}$ modulo $x$.
\end{itemize}
Thus we can compute $\Tr\left(\cN(A^*)
(\xi^2+\tau^2+\Delta)^{-\lambda/2-1}_{[-1]}\right)$ as if all the operators 
involved commute:
\begin{eqnarray*}
II(\lambda)&=&\frac{1}{(2\pi)^2}\int_{S^1\times\rz^2}\Tr\left(\cN(A^*)
(\xi^2+\tau^2+\Delta)^{-\lambda/2-1}_{[-1]}\right) d\tau d\theta d\xi\\
&=&\frac{1}{(2\pi)^2}\int_{S^1\times\rz^2}\Tr\left(
\left(-i\xi+\tau\left[\begin{array}{cc}
1 & 0\\0&-1\end{array}\right]\right)\frac{\left(\frac{\lambda}{2}+1\right)
\left(\frac{\lambda}{2}+2\right)}{2}\right.\\
&&\left.
\left(\frac{4}{i}\tau^2\xi+\frac{2}{i}\tau
\tilde{\nabla}_{\partial_\theta}(\Delta)\right) 
(\xi^2+\tau^2+\Delta)^{-\frac{\lambda}{2}-3}
\right)d\tau d\xi d\theta.
\end{eqnarray*}
We first eliminate the terms which are odd in $\xi$ or $\tau$ 
and thus vanish after integration.
The term containing $\tau^2\left[\begin{array}{cc}
1 & 0\\0&-1\end{array}\right]\tilde{\nabla}_{\partial_\theta}(\Delta)$
is also seen to vanish because the traces on $E^+$ and $E^-$ cancel each other. 
We are left with the term containing $-4\xi^2\tau^2$.
\begin{eqnarray}
II(\lambda)&=&-\frac{4}{(2\pi)^2}\frac{\left(\frac{\lambda}{2}+1\right)
\left(\frac{\lambda}{2}+2\right)}{2}
\nonumber\\&&\int_{S^1\times\rz^2} 
\xi^2\tau^2\Tr(\xi^2+\tau^2+\Delta)^{-\frac{\lambda}{2}-3}
d\tau d\xi d\theta\nonumber\\&=&-\frac{1}{4\pi \lambda} \int_{S^1}
\Tr(D^*D)^{-\frac{\lambda}{2}}d\theta.\label{tii}
\end{eqnarray}
(we integrated by parts in $\tau$ and $\xi$ and then used polar coordinates in 
the plane $(\tau,\xi)$.)

\subsection{The terms of type $III$}
These are the terms coming from the coefficient of $x$
in $q(AA^*)$, i.e.\ the matrix  $x\left[\begin{array}{cc}
i\xi+\tau&-i\tilde{\nabla}_{\partial_\theta}(D^*)\\i\tilde{\nabla}_{\partial_\theta}(D)&
i\xi-\tau \end{array}\right].$ Again, it is impossible to compute these
terms before taking the trace, however
the other two factors involved commute (modulo $x$) so the trace behaves as 
if all operators involved commuted. We get
the following contribution to (\ref{etlm}):
\begin{eqnarray*}
III(\lambda)&=&\left(-\frac{\lambda}{2}-1\right)\frac{1}{(2\pi)^2}
\int_{S^1\times\rz^2}\Tr\left(\left[\begin{array}{cc}
-i\xi+\tau & D^*\\D&-i\xi-\tau\end{array}\right]\right.\\&&
\left[\begin{array}{cc}
\xi^2+\tau^2+D^*D& 0\\0&\xi^2+\tau^2+DD^*
\end{array}\right]^{-\frac{\lambda}{2}-2}\\&&\left.\left[\begin{array}{cc}
i\xi+\tau&-i\tilde{\nabla}_{\partial_\theta}(D^*)\\i\tilde{\nabla}_{\partial_\theta}(D)&
i\xi-\tau \end{array}\right]\right)d\tau d\xi d\theta.
\end{eqnarray*}
The middle matrix is diagonal.
Let us look first at the terms coming from the diagonal entries in 
the first and third matrix. They give
\begin{eqnarray}
III_1(\lambda)&=&\frac{\left(-\frac{\lambda}{2}-1\right)}{(2\pi)^2}
\int_{S^1\times\rz^2} (\xi^2+\tau^2)\Tr
(\xi^2+\tau^2+\Delta)^{-\frac{\lambda}{2}-2}d\tau d\xi d\theta\nonumber\\
&=&-\frac{1}{(2\pi)^2}\int_{S^1\times\rz^2} \Tr
(\xi^2+\tau^2+\Delta)^{-\frac{\lambda}{2}-1}d\tau d\xi d\theta\nonumber\\
&=&-\frac{1}{2\pi
\lambda}\int_{S^1}\Tr(\Delta^{-\frac{\lambda}{2}})d\theta.\label{tiii1}
\end{eqnarray}
Finally let us compute the contribution coming from anti-diagonal 
entries:
\begin{eqnarray}
III_2(\lambda)&=&-\frac{i(\lambda+2)}{(2\pi)^2}
\int_{S^1\times\rz^2}\left[
\Tr\left(D^*(\xi^2+\tau^2+DD^*)^{-\frac{\lambda}{2}-2}
\tilde{\nabla}_{\partial_\theta}(D)\right)\right.\nonumber\\
&&\left.-\Tr\left(D(\xi^2+\tau^2+D^*D)^{-\frac{\lambda}{2}-2}
\tilde{\nabla}_{\partial_\theta}(D^*)\right)\right] d\tau d\xi d\theta\nonumber\\
&=&-\frac{i}{4\pi}\int_{S^1}\left(\Tr\left(D^*(DD^*)^{-\frac{\lambda}{2}-1}
\tilde{\nabla}_{\partial_\theta}(D)\right)\right.\nonumber\\&&
-\left.\Tr\left(D(D^*D)^{-\frac{\lambda}{2}-1}
\tilde{\nabla}_{\partial_\theta}(D^*)\right)\right)
d\theta\nonumber\\
&=&-\frac{i}{2\pi}\int_{S^1}\omega^{BF}(\lambda)d\theta.\label{tiii2}
\end{eqnarray}
The terms (\ref{u1}), (\ref{ti}), (\ref{tii}) and (\ref{tiii1}) cancel so
Proposition \ref{prop8.4} follows from 
(\ref{tiii2}) specialized at $\lambda=0$ and from (\ref{loghol}).
\qed

By Theorem \ref{ldal}
the quantity computed in Proposition \ref{prop8.4} equals half the
adiabatic limit (the limit as $x$ tends to $0$)
of the eta invariant of the family $\delta_x$. 
To complete the proof of the index Theorem \ref{index} we must show the integrability
of the index density. 

\begin{p}\label{intl1}
Under the assumptions of Theorem {\rm \ref{index}}, the local index density 
$\tr(r_1(0)-r_2(0))$ is a smooth multiple of $1/x$. Moreover, 
$$\lim_{\varepsilon\to 0} \int_{X\cap\{x\geq\varepsilon\}} \tr(r_1(0)-r_2(0))$$
exists and gives the $\overline{AS}$ term in the index formula without regularization
with $x^z$.
\end{p}
\dem By Corollary \ref{kdq} we know that 
\begin{equation}\label{trj}
\tr(r_j(0))\sim_{x\to 0}\frac{1}{(2\pi)^2}\left(\int_{\rz^2} 
\tr q(R_j^{-\lambda})d\xi d\tau\right) \frac{d\theta dx}{x^{3}}
\end{equation}
has a Laurent expansion at $x=0$ with a possible singularity of order $3$.
Thus we first want to show that the coefficients of $x^{-3}$, $x^{-2}$ 
in $\tr(r_1(0))-\tr(r_2(0))$
vanish.
We caution the reader that the products in this proof are with respect to the 
rule (\ref{prodsus}). Since $\tr$ is invariant under conjugation by
linear isomorphisms, we can replace the operator $U$ near $x=0$ with 
$$P:=(x^2\partial_x -\frac{x}{2})I_2+\delta_x.$$
We have 
\begin{eqnarray*}
q(PP^*)&=&\xi^2+ix\xi+\frac{x^2}{2}+
q(\delta_x^2)+x\tau\left[\begin{array}{cc}-1&0\\0&1\end{array}\right]\\
q(P^*P)&=&\xi^2+ix\xi+\frac{x^2}{2}+q(\delta_x^2)
-x\tau\left[\begin{array}{cc}-1&0\\0&1\end{array}\right].
\end{eqnarray*}
From this we see that $q(PP^*)^{-\lambda}=q(P^*P)^{-\lambda}$ modulo $x$, and
that $(q(PP^*)^{-\lambda}-q(P^*P)^{-\lambda})_{[-1]}$ is odd in $\tau$ hence vanishes
after integration. Therefore (\ref{trj}) proves the first part of the theorem.

Let us now examine the integral on $\bX$
of the coefficient of $x^{-1}$ in the index density. 
By (\ref{trj}) this is
$$\int_{\rz^2\times S^1} \Tr(q(PP^*)^{-\lambda}-q(P^*P)^{-\lambda})_{[-2]}d\tau d\xi
d\theta.$$
We cannot give an explicit formula for 
$(q(PP^*)^{-\lambda}-q(P^*P)^{-\lambda})_{[-2]}$, but we can do it for its trace. 
We proceed as in the proof of Proposition \ref{prop8.4} to eliminate the terms odd in
$\tau$ or $\xi$. We are left with 
\begin{eqnarray*}
&\int_{\rz^2\times S^1}&
\frac{-\frac{\lambda}{2}\left(-\frac{\lambda}{2}-1\right)}{2}\\
&&\Tr\left(\frac{1}{i}x^2 \left[\begin{array}{cc}-1&0\\0&1\end{array}\right]
\tilde{\nabla}_{\partial_\theta}(\Delta)
(\xi^2+\tau^2+\Delta)^{-\frac{\lambda}{2}-2}\right)d\xi d\tau d\theta
\end{eqnarray*}
(coming from the composition of
$x\tau\left[\begin{array}{cc}-1&0\\0&1\end{array}\right]$ with
$(\xi^2+\tau^2+\Delta)^{-\frac{\lambda}{2}-1}$)
and
\begin{eqnarray*}
&\int_{\rz^2\times S^1}&
\frac{\left(-\frac{\lambda}{2}\right)\left(-\frac{\lambda}{2}-1\right)
\left(-\frac{\lambda}{2}-2\right)}{2}\\
&&\Tr\left(\frac{1}{i}x^2 2\tau^2\left[\begin{array}{cc}-1&0\\0&1\end{array}\right]
\tilde{\nabla}_{\partial_\theta}(\Delta)
(\xi^2+\tau^2+\Delta)^{-\frac{\lambda}{2}-3}\right)d\xi d\tau d\theta
\end{eqnarray*}
(coming from $(\xi^2+\tau^2+\Delta)^{-\frac{\lambda}{2}-1}_{[-1]}$). Integration by
parts with respect to $\tau$ in the second term gives the negative of the first term.
Note that this canceling occurs before the integration in $\theta$. Both terms are
actually $0$ after integration in all variables.
\qed

Note that we proved slightly more than we claimed, namely that the fiberwise 
integral of the index density is smooth as a density in $x,\theta$. It seems 
reasonable to ask if the index
density itself is smooth down to $x=0$ (as in the case of Dirac operators), 
however we were unable to prove or to disprove this fact.

\bibliographystyle{plain}

\begin{thebibliography}{1}

\bibitem{aps1}
M.F.~Atiyah, V.K.~Patodi, and I.M.~Singer,
{\sl Spectral asymmetry and Riemannian geometry. III., }
Math. Proc. Cambridge Philos. Soc. {\bf 79} (1976), 71--79.

\bibitem{bgv}
N.~Berline, E.~Getzler, M.~Vergne,
{\sl Heat Kernels and Dirac Operators, }
Springer-Verlag, Berlin Heidelberg 1991.

\bibitem{bc}
J.-M.~Bismut and J.~Cheeger,
{\sl $\eta$-invariants and their adiabatic limits, }
J. Am. Math. Soc. {\bf 2} (1989), 33--70.

\bibitem{bf2}
J.-M.~Bismut and D.S.~Freed,
 {\sl The analysis of elliptic families II: Dirac operators, eta
  invariants and the holonomy theorem of Witten, }
 Commun. Math. Phys. {\bf 107} (1986), 103--163.

\bibitem{buc}
B.~Bucicovschi,
 {\sl An extension of the work of V. Guillemin on complex powers
and zeta
functions of elliptic pseudodifferential operators, } 
Proc. Amer. Math. Soc. {\bf 127}, No.10 (1999), 3081--3090.


\bibitem{phiho1}
R.~Lauter and S.~Moroianu,
{\sl Homology of pseudodifferential operators on manifolds with fibered
  cusps, } to appear in T. Am. Math. Soc.

\bibitem{de2}
R.~Lauter and S.~Moroianu,
{\sl Homology of pseudo-differential operators on manifolds with fibered
  boundaries, } Journal Reine Angew. Math. {\bf 547} (2002), 207--234.

\bibitem{seronote}
R.~Lauter and S.~Moroianu,
{\sl On the index formula of Melrose and Nistor, }
IMAR preprint 3/2000, Bucharest.

\bibitem{defr}
R.~Lauter and S.~Moroianu,
{\sl Fredholm theory for degenerate pseudodifferential operators on
  manifolds with fibered boundaries, }
Comm. Partial Differential Equations {\bf 26} (2001), 233--283.

\bibitem{maz91}
R.R.~Mazzeo,
{\sl Elliptic theory of differential edge operators {I}, }
Comm. Partial Differential Equations {\bf 16} (1991), 1615--1664.


\bibitem{mame99}
R.R.~Mazzeo and R.B.Melrose,
{\sl Pseudodifferential operators on manifolds with fibered boundaries, }
Asian J. Math. {\bf 2} (1998), 833--866.

\bibitem{mesca}
R.B.~Melrose,
{\sl Spectral and scattering theory for the {L}aplacian on asymptotically
  {E}uclidean space, }
In M.~Ikawa (ed.), Spectral and Scattering Theory, vol. 162
   Lect. Notes Pure Appl. Math., 85--130, New
  York, 1994. Marcel Dekker Inc.
Proc. Taniguchi International Workshop, Sanda,
  November 1992.

\bibitem{meicm}
R.B.~Melrose,
{\sl Pseudodifferential operators, corners and singular limits, }
in Proc. International Congress of Mathematicians,
  Kyoto, Springer-Verlag Berlin - Heidelberg - New York (1990), 217--234.

\bibitem{meleta}
R.B.~Melrose,
 {\sl The eta invariant and families of pseudodifferential operators}.
Math. Res. Letters {\bf 2} (1995), 541--561.

\bibitem{mecom}
R.B.~Melrose,
{\sl Fibrations, compactifications and algebras of pseudodifferential
  operators, }
in L.~H\"ormander and A.~Mellin, editors, {\em Partial Differential
  Equations and Mathematical Physics} (1996), 246--261.

\bibitem{meni96c}
R.B.~Melrose and V.~Nistor, 
{\sl Homology of pseudodifferential operators {I}. {M}anifolds with
  boundary, }
to appear in Amer.\ Math. J., Preprint, May 1996.


\bibitem{ela}
S.~Moroianu,
 {\sl Sur la limite adiabatique des fonctions \^eta et z\^eta, }
Comptes Rendus Math. {\bf 334}, No.2 (2002), 131--134.


\bibitem{etalim}
S.~Moroianu,
{\sl Adiabatic limits of eta and zeta functions, }
math.DG/0204163, to appear in Math. Z.

\bibitem{nysi}
T.M.W.~Nye and M.A.~Singer,
{\sl An ${L}\sp 2$-index theorem for {D}irac operators on ${S}\sp
  1\times{\rz}\sp 3$, } 
J. Funct. Anal. {\bf 177} (2000), 203--218.

\bibitem{boristh}
B.~Vaillant,
{\sl Index- and spectral theory for manifolds with generalized
  fibered cusps, }
PhD thesis, University of Bonn (2001).

\bibitem{wit}
E.~Witten,
{\sl Global gravitational anomalies, }
Commun. Math. Phys. {\bf 100} (1985), 197--229.


\end{thebibliography}

\end{document}